\documentclass[12pt]{report} 
\usepackage{latexsym}
\newcommand{\insp}{\hspace*{1cm}}

\newcommand{\ds}{\displaystyle}

\usepackage{amsfonts}

\usepackage{epsfig}

\newcommand{\ad}{\mbox{ad}}

\newcommand{\reg}{\mbox{\scriptsize{reg}}}

\newcommand{\trdeg}{\mbox{trdeg}}

\newcommand{\End}{\mbox{End}}
\newcommand{\rank}{\mbox{rank}}

\newcommand{\IR}{I\!\!R}

\newcommand{{\CB}}{\cal B}

\newcommand{{\CC}}{\cal C}

 1
 1

\def\picture#1 by #2 (#3){
\vbox to #2{
\hrule width #1 height 0pt depth 0pt
\vfill
\special{picture #3}}}

\begin{document}
\begin{center}
{\bf \Large On Dixmier's fourth problem}\\
\ \\
{\bf \large Alfons I. Ooms}\\
\ \\
{\it Mathematics Department, Hasselt University, Agoralaan, Campus Diepenbeek, 3590 Diepenbeek, Belgium\\
E-mail address: alfons.ooms@uhasselt.be}
\end{center}

\ \\
{\bf Abstract.}\\
Let $\mathfrak{g}$ be a finite dimensional Lie algebra over an algebraically closed field $k$ of characteristic zero.  Denote by $U(\mathfrak{g})$ its enveloping algebra with quotient division ring $D(\mathfrak{g})$.  In 1974, at the end of his book ``Alg\`ebres enveloppantes'', Jacques Dixmier listed 40 open problems, of which the fourth one asked if the center $Z(D(\mathfrak{g}))$ is always a purely transcendental extension of $k$.  We show this is the case if $\mathfrak{g}$ is algebraic whose Poisson semi-center $Sy(\mathfrak{g})$ is a polynomial algebra over $k$.  This can be applied to many biparabolic (seaweed) subalgebras of semi-simple Lie algebras.\\
We also provide a survey of Lie algebras for which Dixmier's problem is known to have a positive answer.  This includes all Lie algebras of dimension at most 8. We prove this is also true for all 9-dimensional algebraic Lie algebras.  Finally, we improve the statement of Theorem 53 of \cite{45}.
\\
\ \\
{\bf \large 1. Introduction}\\
Let $\mathfrak{g}$ be a Lie algebra over an algebraically closed field $k$ of characteristic zero, with basis $x_1,\ldots, x_n$. Let $U(\mathfrak{g})$ be its enveloping algebra with center $Z(U(\mathfrak{g}))$ and semi-center $Sz(U(\mathfrak{g}))$, i.e. the subalgebra of $U(\mathfrak{g})$ 
generated by the semi-invariants of $U(\mathfrak{g})$.  Denote by $D(\mathfrak{g})$ the quotient division ring of $U(\mathfrak{g})$ with center $Z(D(\mathfrak{g}))$.  
In this paper we address the following problem:\\
Is $Z(D(\mathfrak{g}))$ always rational over $k$ ? (Dixmier's fourth problem [13, p.354]).\\
This is true for a wide range of Lie algebras and we are not aware of the existence of a counterexample. See Section 2.\\
In order to simplify things we consider the symmetric algebra $S(\mathfrak{g})$ which we identify with the polynomial algebra $k[x_1,\ldots, x_n]$. This allows us to use MAPLE for the less trivial calculations.\\
We equip $S(\mathfrak{g})$ with its natural Poisson structure.  Its Poisson center $Y(\mathfrak{g})$ coincides with the algebra $S(\mathfrak{g})^{\mathfrak{g}}$ of invariants.
By a celebrated result of Michel Duflo \cite{15}, \cite{16}, \cite{17} there exists an algebra isomorphism between $Z(U(\mathfrak{g}))$ and $Y(\mathfrak{g})$.  
Rentschler and Vergne \cite{50} later extended this to an algebra isomorphism between $Sz(U(\mathfrak{g}))$
and the semi-center $Sz(S(\mathfrak{g}))$, which is usually denoted by $Sy(\mathfrak{g})$.  Furthermore, $Z(D(\mathfrak{g}))$ is isomorphic 
with $R(\mathfrak{g})^{\mathfrak{g}}$, the subfield of invariants of
$R(\mathfrak{g})$, where $R(\mathfrak{g})$ is the quotient field of $S(\mathfrak{g})$.  Therefore Dixmier's problem is equivalent with:
\begin{center}
Is $R(\mathfrak{g})^{\mathfrak{g}}$ always rational over $k$ ?
\end{center}
We now collect some results in order to make this paper as self-contained as possible.  The insiders may very well skip this part. For more information see \cite{45}.\\
\ \\
{\bf 1.1 $i(\mathfrak{g})$, the index of $\mathfrak{g}$}\\
For each $\xi \in \mathfrak{g}^\ast$ we consider its stabilizer
\begin{eqnarray*}
\mathfrak{g}(\xi) = \{x \in \mathfrak{g}\mid \xi ([x,y]) = 0\ \mbox{for all}\ y \in \mathfrak{g}\}
\end{eqnarray*}
The minimal value of $\dim \mathfrak{g}(\xi)$ is called the index of $\mathfrak{g}$ and is denoted by $i(\mathfrak{g})$ [13, 1.11.6; 52, 19.7.3].
Put $c(\mathfrak{g}) = (\dim \mathfrak{g} + i(\mathfrak{g}))/2$.  This integer is sometimes called the magic number of $\mathfrak{g}$.
An element $\xi \in \mathfrak{g}^\ast$ is called regular if $\dim \mathfrak{g}(\xi) = i(\mathfrak{g})$.  
The set $\mathfrak{g}^\ast_{\reg}$ of all regular elements of $\mathfrak{g}^\ast$ is an open dense subset of $\mathfrak{g}^\ast$.\\
We recall from [13, 1.14.13] that
$$i(\mathfrak{g}) = \dim \mathfrak{g} - \rank_{R(\mathfrak{g})} ([x_i, x_j])$$
In particular, $\dim \mathfrak{g} - i(\mathfrak{g})$ is an even number.\\
\ \\
{\bf Theorem 1.1}
$$\mbox{trdeg}_kR(\mathfrak{g})^{\mathfrak{g}}= \mbox{trdeg}_kZ(D(\mathfrak{g})) \leq i(\mathfrak{g})$$
Moreover, equality occurs if one of the following conditions is satisfied:
\begin{itemize}
\item[(1)] $\mathfrak{g}$ is algebraic [38, 39; 50 p.401]
\item[(2)] $\mathfrak{g}$ has no proper semi-invariants in $S(\mathfrak{g})$ (or equivalently in $U(\mathfrak{g})$) [47, Proposition 4.1].
\end{itemize}

{\bf 1.2 Commutative polarizations of $\mathfrak{g}$}\\
Suppose $\mathfrak{g}$ admits a commutative Lie subalgebra $\mathfrak{h}$ such that $\dim \mathfrak{h} = c(\mathfrak{g})$, i.e. $\mathfrak{h}$ is a commutative polarization 
(notation: CP) with respect to any $\xi \in \mathfrak{g}_{\reg}^\ast$ [13, 1.12].\\
These CP's occur frequently in the nilpotent case \cite{43},\cite{44}.  If in addition $\mathfrak{h}$ is an ideal of $\mathfrak{g}$ then we call $\mathfrak{h}$ a CP-ideal (notation: CPI).
If a solvable Lie algebra $\mathfrak{g}$ admits a CP then it also admits a CPI [22, Theorem 4.1]. \\
\ \\
{\bf 1.3 The Poisson algebra $S(\mathfrak{g})$ and its center}\\
The symmetric algebra $S(\mathfrak{g})$, which we identify with $k[x_1,\ldots, x_n]$, has a natural Poisson algebra structure, the Poisson bracket of $f, g \in S(\mathfrak{g})$ given by:
$$\{f,g\} = \sum\limits_{i=1}^n \sum\limits_{j=1}^n [x_i,x_j] \ds\frac{\partial f}{\partial x_i} \ds\frac{\partial g}{\partial x_j}$$
In particular, $S(\mathfrak{g}),\{,\}$ is a Lie algebra for which $\mathfrak{g}$ is a Lie subalgebra since for any two elements $x, y \in \mathfrak{g}$ we have that $\{x,y\} = [x,y]$.  Also, for all $f,g,h \in S(\mathfrak{g})$:
$$\{f, gh\} = \{f,g\}h + g\{f,h\} \insp (\ast)$$
It now easily follows that the center $Y(\mathfrak{g})$ of $S(\mathfrak{g}),\{,\}$ is equal to
$$\{f \in S(\mathfrak{g})\mid \{x,f\} = 0\ \ \forall x \in \mathfrak{g}\}$$
and since $\{x,f\} = \ad~x(f)$ this clearly coincides with $S(\mathfrak{g})^{\mathfrak{g}}$, the subalgebra of invariant polynomials of 
$S(\mathfrak{g})$.\\
The Poisson bracket has a unique extension to the quotient field $R(\mathfrak{g})$ of $S(\mathfrak{g})$ such that $(\ast)$ holds in $R(\mathfrak{g})$.  
It follows that $R(\mathfrak{g}),\{,\}$ is a Lie algebra with center $R(\mathfrak{g})^{\mathfrak{g}}$, the subfield of rational invariants of 
$R(\mathfrak{g})$. $R(\mathfrak{g})$ is called the rational Poisson algebra [54, p. 311].\\
If $\mathfrak{h}$ is a CP of $\mathfrak{g}$ then $R(\mathfrak{h})$ is a maximal Poisson commutative 
subfield of $R(\mathfrak{g})$ [41, Theorem 14].\\
\ \\
{\bf 1.4 The semi-center $Sy(\mathfrak{g})$ of $S(\mathfrak{g})$}\\
Let $\lambda \in \mathfrak{g}^\ast$.  We denote by $S(\mathfrak{g})_{\lambda}$ the set of all $f \in S(\mathfrak{g})$ such that $\ad~x(f) = \lambda(x)f$ for all 
$x \in \mathfrak{g}$.  Any element $f \in S(\mathfrak{g})_{\lambda}$ is said to be a semi-invariant w.r.t. the weight $\lambda$.  We call $f$ a proper semi-invariant 
if $\lambda \neq 0$.  Clearly, $S(\mathfrak{g})_{\lambda} S(\mathfrak{g})_{\mu} \subset S(\mathfrak{g})_{\lambda + \mu}$ for all $\lambda, \mu \in \mathfrak{g}^\ast$.  
Let $f, g \in S(\mathfrak{g})$.  If $fg$ is a nonzero semi-invariant of $S(\mathfrak{g})$, then so are $f$ and $g$. Note that $S(\mathfrak{g})_0 = Y(\mathfrak{g})$.\\
The sum of all $S(\mathfrak{g})_{\lambda}$, $\lambda \in \mathfrak{g}^\ast$, is direct and it is a nontrivial factorial subalgebra $Sy(\mathfrak{g})$ of $S(\mathfrak{g})$ \cite{11}, \cite{35}, \cite{36}. Clearly $Sy(\mathfrak{g})$ is graded by $S(\mathfrak{g})_\lambda$, $\lambda \in \mathfrak{g}^\ast$. Moreover, it is Poisson commutative [47, p. 308].\\
Any nonzero semi-invariant can be written uniquely as a product of irreducible semi-invariants.\\
\ \\
{\bf Proposition 1.2}\\
Suppose $h \in R(\mathfrak{g})$, $h \neq 0$.  Then $h \in R(\mathfrak{g})^{\mathfrak{g}}$ if and only if $h$ can be written as a quotient of two semi-invariants of the same weight \cite{11}.\\
Assume that $\mathfrak{g}$ has no proper semi-invariants (as it is if $[\mathfrak{g},\mathfrak{g}]= \mathfrak{g}$).  Then $R(\mathfrak{g})^{\mathfrak{g}}$ is the quotient field of $S(\mathfrak{g})^{\mathfrak{g}} = Y(\mathfrak{g})$.  In particular,
$$\trdeg_k Y(\mathfrak{g}) = \trdeg_k R(\mathfrak{g})^{\mathfrak{g}} = i(\mathfrak{g})$$
by Theorem 1.1.  Also, $\mathfrak{g}$ is unimodular (i.e. $\mbox{tr}(\ad~x) = 0$ for all $x \in \mathfrak{g}$) by [14, Thm. 1.11] and its proof.\\
\ \\
{\bf 1.5 The canonical truncation $\mathfrak{g}_\Lambda$}\\
The weights of the semi-invariants of $S(\mathfrak{g})$ form an additive semi-group $\Lambda(\mathfrak{g})$, which is not necessarily finitely generated [14, p. 322].\\
Next, we denote by $\mathfrak{g}_\Lambda$ the intersection of $\mbox{ker} \lambda$, $\lambda \in \Lambda(\mathfrak{g})$.  $\mathfrak{g}_\Lambda$ is a characteristic ideal of 
$\mathfrak{g}$ which contains $[\mathfrak{g},\mathfrak{g}]$.  It is called the canonical truncation of $\mathfrak{g}$. In most cases $\mathfrak{g}_{\Lambda}$ has a nicer structure than $\mathfrak{g}$.\\
We now apply the same approach as in [10, pp. 331-334] and [37, pp. 213-214] in order to obtain the following.   
In fact (1), (2), (3) do not require for $k$ to be algebraically closed.  See also \cite{5}, \cite{23}, \cite{25}, \cite{50}.\\
\ \\
{\bf Theorem 1.3}
\begin{itemize}
\item[1.] The centralizer $C(Sy(\mathfrak{g})) = S(\mathfrak{g}_\Lambda)$ and $ C_R(Sy(\mathfrak{g})) = R(\mathfrak{g}_\Lambda)$
\item[2.] $\mathfrak{g}_\Lambda$ has no proper semi-invariants and so $R(\mathfrak{g}_\Lambda)^{\mathfrak{g}_\Lambda}$ is the quotient field of $Y(\mathfrak{g}_\Lambda)$.
Also $\mbox{trdeg}_kY(\mathfrak{g}_\Lambda) = i(\mathfrak{g}_\Lambda)$
\item[3.] $S(\mathfrak{g})^{\mathfrak{g}_\Lambda} = Y(\mathfrak{g}_\Lambda)$ and $R(\mathfrak{g})^{\mathfrak{g}_\Lambda} = R(\mathfrak{g}_\Lambda)^{\mathfrak{g}_\Lambda}$
\item[4.] $c(\mathfrak{g}_\Lambda) = c(\mathfrak{g})$ (use [47, Lemma 3.7] and [43, Proposition 3.2])
\item[5.] $Sy(\mathfrak{g}) \subset Y(\mathfrak{g}_\Lambda) = Sy(\mathfrak{g}_\Lambda)$ and equality occurs if $\mathfrak{g}$ is almost algebraic or if $\mathfrak{g}$ 
is Frobenius (i.e. $i(\mathfrak{g}) = 0$)
\item[6.] Suppose $\mathfrak{h}$ is a CP-ideal of $\mathfrak{g}$. Then $\mathfrak{h} \subset \mathfrak{g}_\Lambda$ [22, p. 141] and $Y(\mathfrak{g}_\Lambda) 
\subset S(\mathfrak{h})$
\end{itemize}

{\bf 1.6 The fundamental semi-invariant $p_\mathfrak{g}$}\\
{\bf Definition} Put $t = \dim \mathfrak{g} - i(\mathfrak{g})$, which is the rank of the structure matrix $B = ([x_i,x_j]) \in M_n(R(\mathfrak{g}))$, 
where $x_1, \ldots, x_n$ is an arbitrary basis of $\mathfrak{g}$.  Assume first that $\mathfrak{g}$ is nonabelian.  Then the greatest common divisor 
$q_{\mathfrak{g}}$ of the $t \times t$ minors in $B$ is a nonzero semi-invariant of $S(\mathfrak{g})$ [10, pp. 336-337].  If $\mathfrak{g}$ is abelian we put 
$q_{\mathfrak{g}} = 1$.  Next, let $p_{\mathfrak{g}}$ be the greatest common divisor of the Pfaffians of the principal $t \times t$ minors in $B$.  
In particular, $\deg p_{\mathfrak{g}} \leq (\dim \mathfrak{g} - i(\mathfrak{g}))/2$.  By [47, Lemma 2.1] $p_{\mathfrak{g}}^2 = q_{\mathfrak{g}}$ up to a 
nonzero scalar multiplier.  We call $p_{\mathfrak{g}}$ the fundamental semi-invariant of $S(\mathfrak{g})$ (instead of $q_{\mathfrak{g}}$ as we did in [47, p. 309)].\\
\ \\
{\bf Lemma 1.4} [30, Lemma 2.3] Let $\mathfrak{g}$ be an algebraic Lie algebra. Then $p_{\mathfrak{g}_\Lambda}$ divides $p_{\mathfrak{g}}$.\\
\ \\
{\bf 1.7 Frobenius Lie algebras}\\
A Lie algebra $\mathfrak{g}$ is called Frobenius if there is a linear functional $\xi \in \mathfrak{g}^\ast$ such that the alternating bilinear form
$B_\xi (x,y) = \xi([x,y]), x, y \in \mathfrak{g}$, is nondegenerate, i.e. $i(\mathfrak{g}) = 0$.  In particular, $R(\mathfrak{g})^{\mathfrak{g}} = k$ by Theorem 1.1.
Frobenius Lie algebras form a large class and they appear naturally in different areas.  For example many parabolic and biparabolic (seaweed) subalgebras of semi-simple Lie algebras are Frobenius 
\cite{7}, \cite{8}, \cite{9}, \cite{18}, \cite{19}, \cite{20}, \cite{21}, \cite{33},  \cite{40}, \cite{49}, including most Borel subalgebras of simple Lie algebras [22, p. 146].  
A Frobenius biparabolic Lie algebra $\mathfrak{g}$ satisfies interesting properties.  
For instance $\mathfrak{g}_\Lambda = [\mathfrak{g},\mathfrak{g}]$ [30, Proposition 7.6].\\
\ \\
We now collect some useful facts on semi-invariants from \cite{10}, \cite{40}.  Let $\mathfrak{g}$ be a Frobenius Lie algebra with basis $x_1,\ldots, x_n$.  
Then $n$ is even and $\mathfrak{g}$ has a trivial center.  The Pfaffian $Pf([x_i, x_j]) \in S(\mathfrak{g})$ is homogeneous of degree $\ds\frac{1}{2} \mbox{dim}\ \mathfrak{g}$ 
and $(Pf([x_i, x_j]))^2 = \mbox{det} ([x_i, x_j]) \neq 0$ by 1.1.  Hence $p_\mathfrak{g} = Pf([x_i, x_j])$.  We put $\Delta(\mathfrak{g}) = \mbox{det} ([x_i,x_j])$ (which is well determined up to
nonzero scalar multipliers).\\
$p_\mathfrak{g}$ is a semi-invariant with weight $\tau$, where $\tau(x) = tr(\ad~x)$, $x \in \mathfrak{g}$.\\
\ \\
{\bf Theorem 1.5} Let $\mathfrak{g}$ be Frobenius.  Decompose $p_\mathfrak{g}$ into a product of irreducible factors:
$$p_\mathfrak{g} = v_1^{m_1} \ldots v_r^{m_r}, \ \ \ m_i \geq 1$$
Then:
\begin{itemize}
	\item[(1)] $v_1,\ldots,v_n$ are the only (up to nonzero scalars) irreducible semi-invariants of $S(\mathfrak{g})$, say with weights $\lambda_1,\ldots, \lambda_r \in \Lambda(\mathfrak{g})$.
	\item[(2)] $Sy(\mathfrak{g}) = k[v_1, \ldots, v_r] = Y(\mathfrak{g}_\Lambda)$, a polynomial algebra over $k$.
	\item[(3)] $r = \dim\mathfrak{g} - \dim\mathfrak{g}_\Lambda = i(\mathfrak{g}_\Lambda)$
	\item[(4)] $\lambda_1, \ldots, \lambda_r$ are linearly independent over $k$. They generate the semi-group $\Lambda(\mathfrak{g})$ and $\mathfrak{g}_\Lambda = \cap 
	\ \mbox{ker} \lambda_i$, $i = 1,\ldots, r$
	\item[(5)] $\mathfrak{h}_i = \mbox{ker} \lambda_i$ is an ideal of $\mathfrak{g}$ of index one and $Y(\mathfrak{h}_i) = k[v_i]$ and $R(\mathfrak{h}_i)^{\mathfrak{h}_i} = k(v_i)$
	\item[(6)] $m_1 \lambda_1 + \ldots + m_r \lambda_r = \tau $ and $m_1 \deg v_1 + \ldots + m_r \deg v_r = \deg p_\mathfrak{g} = \ds\frac{1}{2} \dim \mathfrak{g}$
	\item[(7)] (Joseph [30, 2.2]) Suppose in addition that $\mathfrak{g}$ is algebraic.  \\
	Then $p_{\mathfrak{g}_\Lambda} = v_1^{m_1-1} \ldots v_r^{m_r-1}$. 
\end{itemize}

{\bf 1.8 The Frobenius semi-radical $F(\mathfrak{g})$}\\
Put $F(\mathfrak{g}) = \sum\limits_{\xi \in \mathfrak{g}_{\reg}^{\ast}} \mathfrak{g}(\xi)$.  
This is a characteristic ideal of $\mathfrak{g}$ containing $Z(\mathfrak{g})$ and for which $F(F(\mathfrak{g})) = F(\mathfrak{g})$.  
It can also be characterized as follows: $R(\mathfrak{g})^{\mathfrak{g}} \subset R(F(\mathfrak{g}))$ and if $\mathfrak{g}$ is algebraic then $F(\mathfrak{g})$ is 
the smallest Lie subalgebra of $\mathfrak{g}$ with this property.  Similar results hold in $D(\mathfrak{g})$ [42, Proposition 2.4, Theorem 2.5].
Also, $F(\mathfrak{g}) \subset \mathfrak{g}_\Lambda$.\\
As a special case we have the following:\\
$Y(\mathfrak{g}) \subset S(F(\mathfrak{g}))$ (respectively $Z(U(\mathfrak{g})) \subset U(F(\mathfrak{g})))$ and $F(\mathfrak{g})$ is the smallest Lie 
subalgebra of $\mathfrak{g}$ with this property in case $\mathfrak{g}$ is an algebraic Lie algebra without proper semi-invariants. $F(\mathfrak{g})$ played an important role in \cite{56}.
\begin{center}
	$F(\mathfrak{g}) = 0$ if and only if $\mathfrak{g}$ is Frobenius\\
\end{center}
For this reason $F(\mathfrak{g})$ is called the Frobenius semi-radical of $\mathfrak{g}$.  If  $\mathfrak{g}$ admits a CP  $\mathfrak{h}$ then $F(\mathfrak{g})$ is commutative (since $F(\mathfrak{g}) \subset  \mathfrak{h}$).\\
We call $\mathfrak{g}$ quasi-quadratic if $F(\mathfrak{g}) = \mathfrak{g}$.  Such a Lie algebra does not possess any proper semi-invariant.  For example any quadratic Lie algebra is quasi-quadratic.\\
\ \\
{\bf 1.9 Necessary conditions for polynomiality}\\
\ \\
{\bf Theorem 1.6} [47, Theorem 1.1] Let $\mathfrak{g}$ be a Lie algebra for which the semi-center $Sy(\mathfrak{g})$ is freely generated by homogeneous elements 
$f_1, \ldots, f_r$.\\
Then
$$\sum\limits_{i=1}^r \deg f_i \leq c(\mathfrak{g})$$

{\bf Definition.} A Lie algebra $\mathfrak{g}$ is called coregular if $Y(\mathfrak{g})$ is a polynomial algebra over $k$.\\
\ \\
{\bf Theorem 1.7} [47, Proposition 1.4]. For an extension see [33, Theorem 2.2].\\
Assume that $\mathfrak{g}$ has no proper semi-invariants and that $Y(\mathfrak{g})$ is freely generated by homogeneous elements $f_1,\ldots, f_r$.  Then
$$\sum\limits_{i=1}^r \deg f_i = c(\mathfrak{g}) - \deg p_{\mathfrak{g}}$$

{\bf Corollary 1.8} Assume that $\mathfrak{g}$ has no proper semi-invariants and that $Y(\mathfrak{g})$ is freely generated by homogeneous elements $f_1,\ldots, f_r$.  Then
$$3i(\mathfrak{g}) + 2\deg p_{\mathfrak{g}} \leq \dim\mathfrak{g} + 2 \dim Z(\mathfrak{g})$$
Moreover, equality occurs if and only if $\deg f_i \leq 2$, $i:1,\ldots, r$.\\
\ \\
{\bf 1.10 A sufficient condition for polynomiality}\\
\ \\
{\bf Theorem 1.9} [33, 5.7: 51]\\
Assume that $f_1, \ldots, f_r \in Y(\mathfrak{g})$, $r = i(\mathfrak{g})$, are algebraically independent homogeneous invariants such that
$$\sum\limits_{i=1}^r \deg f_i \leq c(\mathfrak{g})) - \deg p_\mathfrak{g}$$
Then, equality holds and $Y(\mathfrak{g}) = k[f_1,\ldots, f_r]$.  In particular, $\trdeg_kY(\mathfrak{g}) = i(\mathfrak{g})$.\\
\ \\
See [43, 45, Theorem 30] for an alternative method.\\
Next we recall [22, Proposition 4.6]:\\
{\bf Proposition 1.10}\\
Let $\mathfrak{g}$ be a Lie algebra over $k$ and $W$ a $\mathfrak{g}$-module with $\dim \mathfrak{g} \leq \dim W$. For each $f \in W^\ast$ we put $ \mathfrak{g}(f) = \{ x \in  \mathfrak{g}\mid f(xw) = 0\ \mbox{for all}\ w\in W\}$ the stabilizer of $f$.  Consider the semi-direct product $L =  \mathfrak{g} \oplus W$ in which $[x,w] = xw$, $x \in  \mathfrak{g}$, $w \in W$, and in which $W$ is an abelian ideal.  Then the following are equivalent:
\begin{itemize}
\item[(1)] $ \mathfrak{g}(f) = 0$ for some $f \in W^\ast$
\item[(2)] $i(L) = \dim W - \dim  \mathfrak{g}$
\item[(3)] $W$ is a commutative polarization of $L$
\item[(4)] $R(W)$ is a maximal Poission commutative subfield of $R(L)$.
\end{itemize}
Moreover, if these conditions are satisfied then $W$ is a faithful $ \mathfrak{g}$-module.  Also,
\begin{center}
$Y(L) = S(W)^{\mathfrak{g}}$ and $ R(L)^L = R(W)^{\mathfrak{g}}$
\end{center}

{\bf Proposition 1.11}\\
Let $L$ be a Lie algebra over $k$ having an element $u \in L$ such that the centralizer $M = C(u)$ has codimension one in $L$.  Then $F(L) \subset M$.\\
\ \\
{\bf Proof.} First we obtain that $i(M) = i(L) + 1$ from [22, Proposition 1.9(i]. But this implies that $F(L) \subset M$ by [22, Proposition 1.6(4)].\hfill $\square$\\
\ \\
{\bf \large 2. Some known positive results about Dixmier's 4th problem}\\
\ \\
{\bf Theorem 2.1}  The following finite dimensional Lie algebras over $k$ satisfy Dixmier's 4th problem.
\begin{itemize}
	\item[(1)] Any solvable Lie algebra (Dixmier, [13, 4.4.8]). This remains true if $k$ is replaced by $\IR$ \cite{3}.
	\item[(2)] Any coregular Lie algebra $L$ without proper semi-invariants (Indeed, $Y(L)$ is polynomial and $R(L)^L$ is the quotient field of $Y(L)$ by Proposition 1.2).
	\item[(3)] Any semi-simple Lie algebra $L$ (This is a special case of (2) as $Y(L)$ is polynomial by a theorem of Chevalley [13, 7.3.8], while $L$ has no proper semi-invariants as $[L,L] = L$).
	\item[(4)] Any Frobenius Lie algebra $L$ (as $R(L)^L = k$) and its truncation $L_\Lambda$ (this is also a special case of (2) since $L_\Lambda$ has no proper semi-invariants and $Y(L_\Lambda) = Sy(L)$ by (5) of Theorem 1.3 and $Sy(L)$ is polynomial by (2) of Theorem 1.5).
	\item[(5)] Any biparabolic subalgebra $P$ \cite{48} and its truncation $P_\Lambda$ \cite{28}, \cite{29} of a simple Lie algebra of type $A$ or $C$.
	\item[(6)] Any maximal parabolic subalgebra $P$ of a simple Lie algebra of type $B$, $D$ or $E_6$ \cite{26}.
	\item[(7)] Any algebraic Lie algebra $L$ with adjoint group $G$ having an affine slice $V \subset L^\ast$ for the coadjoint action of $L$ \cite{53}.\\
	Such an affine slice exists for the coadjoint action for certain truncated \cite{31},\cite{32} and non truncated \cite{53} biparabolic subalgebras of a semi-simple Lie algebra.
	\item[(8)] Any square integrable Lie algebra $L$ (i.e. $i(L) = \dim Z(L)$ and thus $F(L) = Z(L)$ which implies that $R(L)^L = R(Z(L))$ [42, Corollary 2.6)].  For an extension of this see [45, Theorem 37].
	\item[(9)] The semi-direct product $L = sl(2)\oplus W_m$, where $W_m$ is the $(m+1)$-irreducible $sl(2)$-module [45, Proposition 62].
	\item[(10)] Any Lie algebra $L$ with $\dim L \leq 8$. [45, 46, Proposition 63].  See also Section~4.
\end{itemize}

The following is equivalent to [45, Theorem 66].  It shows how Dixmier's 4th problem for an algebraic Lie algebra $L$ can be reduced to that of its canonical truncation $L_\Lambda$.\\
\ \\
{\bf Theorem 2.2} Let $L$ be algebraic such that the field $R(L_\Lambda)^{L_\Lambda}$ is freely generated by some semi-invariants $u_1,\ldots, u_s$ of $S(L)$.  Then $R(L)^L$ is rational over $k$.\\
\ \\
{\bf \large 3. More recent results}\\
\ \\
The following was provided by Michel Van den Bergh.\\
{\bf Proposition 3.1}\\
Assume $A$ is a polynomial algebra over $k$ with its natural grading
\begin{center}
	$A = k \oplus A_1 \oplus A_2 \oplus \ldots$ \ \ and let \ \ $M = A_1 \oplus A_2 \oplus \ldots$
\end{center}
be its augmentation ideal.
\begin{itemize}
\item[(1)] Let $u_1,\ldots, u_n$ be a homogeneous lift of a basis of the vector space $M/M^2$.  Then $u_1, \ldots, u_n$ freely generate $A$.
\item[(2)] Suppose $A$ admits an additional grading (as it is if $A$ is the semi-center of a Lie algebra).  Then with respect to this, $u_1,\ldots, u_n$ can be taken to be homogeneous (apply (1) to a homogeneous basis of $M/M^2$).
\end{itemize}

{\bf Corollary 3.2}\\
Let $L$ be a Lie algebra such that its semi-center $Sy(L)$ is polynomial over $k$.  Then $Sy(L)$ is freely generated by some semi-invariants $u_1,\ldots, u_s$ of $S(L)$.\\
\ \\
{\bf Theorem 3.3}\\
Let $L$ be an algebraic Lie algebra such that its Poisson semi-center $Sy(L)$ is polynomial over $k$.  Then $R(L)^L$ is rational over $k$.\\
The converse of this result is not true (see examples $L_{9.10}$ and $L_{9.11}$ of Section 4).\\
\ \\
{\bf Proof.} Since $L$ is algebraic $Y(L_\Lambda) = Sy(L)$ (by (5) of Theorem 1.3), which is polynomial over $k$ by assumption.  So, by the Corollary $Y(L_\Lambda)$ is freely generated by some semi-invariants of $S(L)$, say $u_1,\ldots, u_s$.  Hence $u_1,\ldots, u_s$ freely generate the quotient field of $Y(L_\Lambda)$, which is $R(L_\Lambda)^{L_\Lambda}$ by (2) of Theorem 1.3. Using Theorem 2.2 we may conclude that $R(L)^L$ is rational over $k$. \hfill $\square$\\
\ \\
{\bf Corollary 3.4}\\
In recent years many (bi)parabolic subalgebras of certain semi-simple Lie algebreas have been shown to possess a polynomial Poisson semi-center \cite{23},\cite{24},\cite{25},\cite{26},\cite{28},\cite{29}, \cite{57}.  By the previous Theorem these (bi)parabolic subalgebras satisfy Dixmier's 4th problem.\\
\ \\
{\bf Remark 3.5} In 2007 Oksana Yakimova \cite{55} discovered a counterexample to the polynomiality of the Poisson semi-center of a (bi)parabolic subalgebra $L$, namely when $L$ is a certain maximal parabolic subalgebra of a simple Lie algebra of type $E_8$. At the same time this example is also a counterexample to a conjecture by Premet \cite{55}.  Perhaps it is also a counterexample for Dixmier's 4th problem.\\
\ \\
{\bf Proposition 3.6}\\
Let $L = sl(2) \oplus W$ be the semi-direct product of $sl(2)$ with a $sl(2)$-module $W$ with $\dim W \geq 3$ ($= \dim sl(2)$).  Suppose that $i(L) = \dim W -3$.  Then $W$ is a faithful $sl(2)$-module and a $CP$-ideal of $L$ and $L$ satisfies Dixmier's 4th problem.\\
\ \\
{\bf Proof.} By Proposition 1.10 the first part follows at once and also
$$R(L)^L = R(W)^{sl(2)} = R(W)^{SL(2)}$$
where the latter is rational over $k$ by a result of Katsylo \cite{34}.  See also \cite{4}.\hfill $\square$\\
\ \\
{\bf Proposition 3.7} Assume $L = L_1 \oplus L_2$ is a direct product.  Let $\lambda \in L^\ast$ and put $\lambda_1 = \lambda \mid_{L_1} \in L_1^\ast$ and $\lambda_2 = \lambda \mid_{L_2} \in L_2^\ast$.  Then
\begin{itemize}
\item[(1)] $S(L)_\lambda = S(L_1)_{\lambda_1} S(L_2)_{\lambda_2}$
\item[(3)] $Y(L) = Y(L_1) Y(L_2)$
\item[(3)] $Sy(L) = Sy(L_1) Sy(L_2)$
\end{itemize}

{\bf Proof.} (1) The inclusion $\supset$ is easy to verify.  On the other hand take a nonzero $u \in S(L)_\lambda$.  Let $x_1,\ldots, x_r, y_1, \ldots, y_s$ be a basis of $L$ such that $x_1,\ldots, x_r$ (respectively $y_1,\ldots, y_s)$ is a basis of $L_1$ (resp. of $L_2$).\\
For $m = (m_1,\ldots, m_r)$ and $n = (n_1,\ldots, n_s)$ we denote $x^m = x_1^{m_1} \ldots x_r^{m_r}$ and $y^n = y_1^{n_1} \ldots y_s^{n_s}$.\\
Then $u$ has the following unique decomposition:
\begin{center}
$u = \sum\limits_{m,n} a_{mn} x^my^n = \sum\limits_n p_ny^n$ where $p_n = \sum_m a_{mn} x^m \in S(L_1)$
\end{center}
for some $a_{mn} \in k$.  Next take any $x \in L_1$.  Then $\{x,u\} = \lambda (x)u = \lambda_1(x)u$.  Since $\{x,y^n\} = 0$ we get
$$\sum\limits_n \{x,p_n\} y^n = \{x,u\} = \lambda_1(x)u = \sum\limits_n (\lambda_1(x) p_n)y^n$$
for all $n$.  Therefore $\{x, p_n\} = \lambda_1(x) p_n$ for all $x \in L_1$ and so $p_n \in S(L_1)_{\lambda_1}$ for all $n$.\\
Let $u_1 = p_{n_1},\ldots, u_t = p_{n_t}$ be a $k$-basis of the vector space over $k$ generated by all $p_n$'s.  Then $u = \sum\limits_n p_n y^n $ can be written as $u = \sum\limits_{i=1}^t u_iv_i$ for some $v_i \in k[y_1,\ldots, y_s] = S(L_2)$.\\
\ \\
{\bf Claim.} $u_1,\ldots, u_t$ are also linearly independent over $k[y_1,\ldots, y_s]$.\\
Indeed, suppose $\sum\limits_{i=1}^t u_i w_i = 0$ for some $w_i \in k [y_1,\ldots, y_s]$. \\
Now take any $a_1,\ldots, a_s \in k$.  Then
$$\sum\limits_{i=1}^t u_i w_i (a_1, \ldots, a_s) = 0$$
Since $u_1,\ldots, u_t$ are linearly independent over $k$ and $w_i(a_1,\ldots, a_s) \in k$ it follows that $w_i (a_1,\ldots, a_s) = 0$ for all $a_1,\ldots, a_s \in k$.  Hence $w_i = 0$ for all $i$, establishing the claim.\\
Next take any $y \in L_2$.  Then $\{y,u\} = \lambda(y) u = \lambda_2(y) u$.  As $\{y, u_i\} = 0$, we get
$$\sum\limits_{i=1}^t u_i\{y,v_i\} = \{y,u\} = \lambda_2(y) u = \sum\limits_{i=1}^t u_i (\lambda_2(y) v_i)$$
By the claim we may conclude that for all $i = 1,\ldots, t$ $\{y, v_i\} = \lambda_2(y) v_i$, i.e. $v_i \in S(L_2)_{\lambda_2}$.\\
Finally, $u = \sum\limits_{i=1}^t u_i v_i \in S(L_1)_{\lambda_1} S(L_2)_{\lambda_2}$.\\
\ \\
(2) (take $\lambda = 0$) and (3) follow easily from (1).\hfill $\square$\\
\ \\
{\bf Corollary 3.8.}\\
Assume that $L = L_1 \oplus L_2$ is a direct product and that $R(L_1)^{L_1}$ and $R(L_2)^{L_2}$ are both rational over $k$.  Then the same holds for $R(L)^L$.\\
\ \\
{\bf Proof.}
\begin{itemize}
\item[(1)] First we show that $R(L)^L$ is the field generated by $R(L_1)^{L_1}$ and $R(L_2)^{L_2}$.  It is clear that $R(L)^L$ contains $R(L_1)^{L_1}$ and $R(L_2)^{L_2}$.  Now take a nonzero $z \in R(L)^L$.  Then by Proposition 1.2 $z = uv^{-1}$ where $u, v \in S(L)_\lambda \backslash \{0\}$ for some weight $\lambda \in L^\ast$.  From (1) of the previous Proposition we have $S(L)_\lambda = S(L_1)_{\lambda_1} S(L_2)_{\lambda_2}$ where $\lambda_1 = \lambda \mid_{L_1}$ and $\lambda_2 = \lambda \mid_{L_2}$. Hence,
\begin{center}
$u = \sum\limits_i u_i p_i$ where $u_i \in S(L_1)_{\lambda_1}, p_i \in S(L_2)_{\lambda_2}$ are nonzero
\end{center}
Similarly,
\begin{center}
$v = \sum\limits_j v_j q_j$ where $v_j \in S(L_1)_{\lambda_1}, q_j \in S(L_2)_{\lambda_2}$ are nonzero
\end{center}
\begin{eqnarray*}
z &=& uv^{-1} = \left(\sum\limits_i u_i p_i\right) \left(\sum\limits_j v_j q_j\right)^{-1} = 
\sum\limits_i \left[u_i p_i \left(\sum\limits_j v_j q_j\right)^{-1}\right]\\
&=& \sum\limits_i \left[u_i^{-1}p_i^{-1} \left(\sum\limits_j v_j q_j\right)\right]^{-1} = \sum\limits_i \left[\sum\limits_j \left(u_i^{-1}v_j\right)\left(p_i^{-1} q_j\right)\right]^{-1}
\end{eqnarray*}
where $u_i^{-1} v_j \in R(L_1)^{L_1}$ and $p_i^{-1}q_j \in R(L_2)^{L_2}$ by Proposition 1.2.
\item[(2)] Choose a basis $x_1,\ldots, x_m$, $x_{m+1},\ldots, x_n$ of $L$ such that $x_1,\ldots,x_m$ is a basis of $L_1$ and $x_{m+1},\ldots, x_n$ is a basis of $L_2$.\\
By assumption, there are $z_1, \ldots, z_r$; $z_{r+1},\ldots, z_s$ such that both $R(L_1)^{L_1} = k(z_1,\ldots, z_r)$ and $R(L_2)^{L_2} = k(z_{r+1},\ldots, z_s)$ are purely transcendental extension of $k$.  Using (1) we have that
$$R(L)^L = k(z_1,\ldots, z_r, z_{r+1},\ldots, z_s)$$
So, it remains to show that $z_1,\ldots, z_r, z_{r+1}, \ldots , z_s$ are algebraically independent over $k$, which is equivalent with $\mbox{rank}(J) = s$, where\\
$J = \left(\frac{\partial z_j}{\partial x_i}\right)$, $i= 1,\ldots, n$; $j = 1,\ldots, s$, the Jacobian matrix.\\
\ \\
It is easy to see that
$$J = \left(\begin{array}{cc}
J_1 & O\\O & J_2\end{array}\right)$$
where\\
$J_1 = \left(\frac{\partial z_j}{\partial x_i}\right)$, $i = 1,\ldots, m$; $j = 1,\ldots, r$\\
\ \\
$J_2 = \left(\frac{\partial z_j}{\partial x_i}\right)$, $i = m+1,\ldots, n$; $j = r+1,\ldots, s$\\
\ \\
Then $\rank(J) = \rank (J_1) + \rank (J_2) = r + (s-r) = s$.\hfill $\square$
\end{itemize}
\ \\
{\bf \large 4. Dixmier's 4th problem for Lie algebras of dimension at most 9}\\
First we recall [46, Theorem 53]\\
\ \\
{\bf Theorem 4.1} Let $L$ be a nonsolvable, indecomposable Lie algebra with $\dim L \leq 8$.\\
Then
\begin{itemize}
\item[(1)] $Y(L)$ and $Sy(L)$ are polynomial over $k$
\item[(2)] $R(L)^L$ is rational over $k$
\end{itemize}

{\bf Remark 4.2} The following example shows that the condition that $L$ is indecomposable cannot be removed for part (1) (which we erroneously did in [45, Theorem 53].  This has been corrected in \cite{46}).\\
\ \\
{\bf Example 4.3} Consider the direct product $L = L_1 \oplus L_2$, where $L_1 = sl(2)$ with basis $h, x, y$ with nonzero brackets
$$[h,x] = 2x, [h,y] = -2y, [x,y] = h$$
$L_1$ is coregular as $Y(L_1) = k[c]$ where $c = h^2 + 4xy$, the Casimir element of $L_1$.\\
$L_2$ is the solvable Lie algebra with basis $x_1, x_2, x_3, x_4$ and nonzero brackets
$$[x_1, x_2] = x_2, [x_1, x_3] = x_3, [x_1,x_4] = -2x_4$$
Then $L_2$ is not coregular [45, Example 23].\\
In fact, $Y(L_2) = k [f_1, f_2, f_3]$ where
\begin{center}
$f_1 = x_2^2x_4$, $f_2 = x_3^2x_4$, $f_3 = x_2x_3x_4$ with $f_1f_2 = f_3^2$
\end{center}
Clearly, $L$ is a 7-dimensional nonsolvable Lie algebra.  However it is not coregular since
$$Y(L) = Y(L_1) Y(L_2) = k [c, f_1, f_2, f_3]$$
which is not polynomial. \hfill $\square$\\
\ \\
{\bf Corollary 4.4}\\
Let $L$ be a Lie algebra with $\dim L \leq 8$.  Then $R(L)^L$ is rational over $k$.\\
\ \\
{\bf Proof.}\\
In view of Corollary 3.8 we may assume that $L$ is indecomposable.  Then $R(L)^L$ is rational over $k$:
\begin{itemize}
\item[$\bullet$] by (1) of Theorem 2.1 if $L$ is solvable
\item[$\bullet$] by (2) of Theorem 4.1 if $L$ is nonsolvable.
\end{itemize}
\hfill $\square$\\
\ \\
The following can be derived from [45, Lemma 64] and its proof.\\
\ \\
{\bf Lemma 4.5}\\
Let $L$ be an algebraic Lie algebra with radical $R_1$.  Then there exists a torus $T \subset R_1$ such that $L = L_\Lambda \oplus T$ with $[S,T] = 0$ for a certain Levi subalgebra $S$ of $L$.\\
\ \\
{\bf Theorem 4.6}\\
Let $L$ be a 9-dimensional algebraic Lie algebra.  Then $L$ satisfies Dixmier's 4th problem.\\
\ \\
{\bf Proof.} We may assume that $L$ is indecomposable (by Corollary 3.8) and nonsolvable (by (1) of Theorem 2.1).  We distinguish two cases:\\
\ \\
{\bf (1) \boldmath{$L_\lambda \neq L$}}\\
In particular $\dim L_\Lambda \leq 8$. Since $L$ is algebraic so is $L_\Lambda$ [10, Proposition 1.14] and $Y(L_\Lambda) = Sy(L)$ by (5) of Theorem 1.3.  Because of Theorem 3.3 it suffices to show that $Y(L_\Lambda)$ is polynomial. Let $R_1$ be the radical of $L$.  Invoking Lemma 4.5 there exists a torus $T \subset R_1$ such that $L = L_\Lambda \oplus T$ with $[S,T] = 0$ for a certain Levi subalgebra $S$ of $L$.  Then $S = [S,S] \subset [L,L] \subset L_\Lambda$.  From [6, Corollaire 4 p.91], as $L_\Lambda$ is an ideal of $L$, we know that $R = R_1 \cap L_\Lambda$ is the radical of $L_\Lambda$ and $S = S \cap L_\Lambda$ is a Levi subalgebra of $L_\Lambda$  In particular, 
$$L_\Lambda = S \oplus R$$

{\bf Claim:} $[S,R] \neq 0$ and $S = sl(2)$\\
Suppose $[S,R] = 0$.  Consider $L = L_\Lambda \oplus T = S \oplus (R \oplus T)$ where $R \oplus T$ is the radical $R_1$ of $L$.  Notice that 
$$[S,R \oplus T] = [S,R] + [S,T] = 0$$
Hence $L$ would be decomposable.  Contradiction. $[S,R] \neq 0$ forces $R \neq 0$, i.e. $L_\Lambda$ is not semi-simple.  In particular, $S \neq sl(3)$ (otherwise $L_\Lambda = sl(3)$ as $\dim sl(3) = 8$).\\
Now, suppose $S = S_1 \oplus S_2$ is a direct product of 2 copies of $sl(2)$.  In particular, $\dim S = 6$. $[S,R] \neq 0$ implies that $\dim R = 2$ and
$$\varphi : S \rightarrow \End R, x \mapsto \ad x\mid_R$$
is a nonzero Lie algebra homomorphism, which is not injective (since $\dim S = 6 > 4 = \dim \End R$).\\
Consequently, $\ker \varphi$ is a nontrivial ideal of $S$, say $S_1$, and so $[S_1, R] = 0$, while $[S_2, R] \neq 0$.  Hence, $[S_1, S_2 \oplus R \oplus T] = 0$ which means that $L = S_1\oplus (S_2 \oplus R \oplus T)$ is decomposable.  Contradiction.\\
Therefore $S = sl(2)$, settling the claim.  We may conclude that $L_\Lambda$ is a member of the list of \cite {2}, which we used in order to prove that all algebraic Lie algebras of dimension at most 8, satisfy the Gelfand-Kirillov conjecture \cite{27}.  As a consequence $L_\Lambda$ is either indecomposable and nonsolvable (and thus coregular by Theorem 4.1) or it is a direct product $L_\Lambda = M \oplus N$, where $M$ is indecomposable and nonsolvable (and thus coregular) and $N$ is solvable of one of the following types:
\begin{itemize}
\item[(i)] $N$ is abelian with $\dim N \leq 3$
\item[(ii)] $N$ is the 2-dimensional nonabelian Lie algebra
\item[(iii)] $N$ is the 3-dimensional Heisenberg Lie algebra
\item[(iv)] $N$ is the 3-dimensional Lie algebra with basis $x,y,z$ with brackets:
 \begin{center}
 $[x,y] = y$, $[x,z] = \alpha z$ with $\alpha \in \Bbb Q$
 \end{center}
 \end{itemize}
 For each of the above $N$ is coregular, i.e. $Y(N)$ is polynomial.  Hence the same is true of $Y(L_\Lambda) = Y(M) Y(N)$.\\
 \ \\
{\bf (2) \boldmath{$L_\Lambda = L$}}\\
We will proceed case by case using the classification of the 9-dimensional nonsolvable Lie algebras provided to us by Boris Komrakov.  From this list we select all algebraic, indecomposable Lie algebras $L$ for which $L_\Lambda = L$ (i.e. $L$ has no proper semi-invariants).  In doing so, there are only 11 cases that remain to be examined.\\
\ \\
{\bf Notation:} $x_1, x_2, x_3, x_4,\ldots, x_9$ will be a basis of $L$ such that $x_1, x_2, x_3$ is the standard basis of $sl(2)$ (i.e. $[x_1,x_2] = 2x_2$, $[x_1, x_3] = -2x_3$ and $[x_2, x_3] = x_1$) and $x_4,\ldots, x_9$ is a basis of the radical $W$ of $L$.\\
For each $L$ the structure matrix $M = ([x_i, x_j])$ of $L$ will be given explicitly.\\
We now distinguish two subcases.\\
\ \\
{\bf (2a) \boldmath {$L$} is coregular}\\
In this situation we have at once that $R(L)^L$ is rational over $k$ by (2) of Theorem 2.1 since $L$ has no proper semi-invariants.\\
In each of the following cases our approach is more or less the same.  For this reason we will only give the explicit description of the solution for the first Lie algebra.\\
The following 3 Lie algebras are the only ones in our list for which $[L, L] \neq L$.\\
\ \\
{\bf (1) \boldmath{$L_{9,1}$}}
\begin{eqnarray*}
M := \left[\begin{array}{ccccccccc}
0 &2x_2 &-2x_3 &0 &0 &x_6 &-x_7 &0 &0\\
-2x_2 &0 &x_1 &0 &0 &0 &x_6 &0 &0\\
2x_3 &-x_1 &0 &0 &0 &x_7 &0 &0 &0\\
0 &0 &0 &0 &0 &0 &0 &0 &0\\
0 &0 &0 &0 &0 &0 &0 &0 &x_4\\
-x_6 &0 &-x_7 &0 &0 &0 &x_4 &0 &0\\
x_7 &-x_6 &0 &0 &0 &-x_4 &0 &0 &0\\
0 &0 &0 &0 &0 &0 &0 &0 &x_5\\
0 &0 &0 &0 &-x_4 &0 &0 &-x_5 &0
\end{array}\right]
\end{eqnarray*}
$i(L) = \dim L - \rank (M) = 9-6 = 3$ (see 1.1).\\
$c(L) = (\dim L + i(L))/2 = 6$.  Using MAPLE we get $p_L = 1$.\\
Next we want to determine $F(L)$.  Take $\xi = x_1^\ast + x_4^\ast \in L^\ast$.  Then $L(\xi) = \langle x_1, x_4, x_8\rangle$.\\
Clearly $\dim L(\xi) = 3 = i(L)$ and so $\xi \in L_{\reg}^\ast$.  Therefore $\langle x_1, x_4, x_8\rangle = L(\xi) \subset F(L)$.  Since $F(L)$ is an ideal of $L$ (1.8) we obtain:
$$\langle x_2, x_3, x_6, x_7\rangle = [L, x_1] \subset F(L), \langle x_5\rangle = [L, x_8] \subset F(L)$$
Hence $\langle x_1,\ldots, x_8\rangle \subset F(L)$.  On the other hand, the centralizer $C(x_8) = \langle x_1,\ldots, x_8\rangle$ is of codimension one in $L$ which implies that $F(L) \subset C(x_8)$ by Proposition 1.11 and therefore $F(L) = \langle x_1,\ldots, x_8\rangle$, which is not commutative and so $L$ does not contain any CP's (1.8).\\
Moreover, $\langle x_1,\ldots, x_8\rangle = F(L) \subset L_\Lambda$ (1.8).  Also $x_9 \in L_\Lambda$ since $\ad x_9$ is nilpotent (in fact $(\ad x_9)^3 = 0$).\\
Consequently, $L_\Lambda = L$.\\
Next, using MAPLE, we' obtain the following $3(=i(L))$ homogeneous, algebraically independent (over $k$) elements of $Y(L)$:\\
$x_4$, $f=2x_4x_8 - x_5^2$ and\\
$g = x_1^2x_4 + 2x_1x_6x_7 + 4x_2x_3x_4 + 2x_2x_7^2 - 2x_3x_6^2$\\
Now we observe that:
$$\deg x_4 + \deg f + \deg g = 6 = c(L) - \deg p_L$$
Invoking Theorem 1.9 we may conclude that:
\begin{center}
$Y(L) = k[x_4, f, g]$ which is polynomial over $k$
\end{center}
and $R(L)^L = k(x_4, f, g)$ which is rational over $k$.\\
\ \\
{\bf (2) \boldmath{$L_{9,2}$}}
\begin{eqnarray*}
M := \left[\begin{array}{ccccccccc}
0 &2x_2 &-2x_3 &0 &x_5 &-x_6 &0 &0 &0\\
-2x_2 &0 &x_1 &0 &0 &x_5 &0 &0 &0\\
2x_3 &-x_1 &0 &0 &x_6 &0 &0 &0 &0\\
0 &0 &0 &0 &0 &0 &0 &0 &0\\
-x_5 &0 &-x_6 &0 &0 &x_4 &0 &0 &0\\
x_6 &-x_5 &0 &0 &-x_4 &0 &0 &0 &0\\
0 &0 &0 &0 &0 &0 &0 &x_4 &-x_7\\
0 &0 &0 &0 &0 &0 &-x_4 &0 &x_8\\
0 &0 &0 &0 &0 &0 &x_7 &-x_8 &0
\end{array}\right]
\end{eqnarray*}
$i(L) = 3$, $c(L) = 6$, $p_L = 1$, $F(L) = L$ (so $L$ is quasi quadratic, with no CP's). This implies that $L_\Lambda = L$.\\
$Y(L) = k[x_4, f, g]$ which is polynomial over $k$, where $f = x_4 x_9 + x_7x_8$,\\
$g = x_1^2x_4 + 2x_1x_5x_6 + 4x_2x_3x_4 + 2x_2x_6^2 - 2x_3x_5^2$\\
$R(L)^L = k(x_4, f, g)$ which is rational over $k$.\\
\ \\
{\bf (3) \boldmath{$L_{9,3}$}}
\begin{eqnarray*}
M := \left[\begin{array}{ccccccccc}
0 &2x_2 &-2x_3 &0 &x_5 &-x_6 &x_7 &-x_8 &0\\
-2x_2 &0 &x_1 &0 &0 &x_5 &0 &x_7 &0\\
2x_3 &-x_1 &0 &0 &x_6 &0 &x_8 &0 &0\\
0 &0 &0 &0 &0 &0 &0 &0 &0\\
-x_5 &0 &-x_6 &0 &0 &0 &0 &0 &0\\
x_6 &-x_5 &0 &0 &0 &0 &0 &0 &0\\
-x_7 &0 &-x_8 &0 &0 &0 &0 &x_4 &x_5\\
x_8 &-x_7 &0 &0 &0 &0 &-x_4 &0 &x_6\\
0 &0 &0 &0 &0 &0 &-x_5 &-x_6 &0
\end{array}\right]
\end{eqnarray*}
$i(L) = 3$, $c(L) = 6$, $p_L = 1$, $F(L) = L$ (so $L$ is quasi quadratic with no CP's), $[ L,L] = L$, $L_\Lambda = L$.\\
$Y(L) = k[x_4, f,g]$ which is polynomial over $k$, where\\
$f = x_4x_9 - x_5x_8 + x_6x_7$ and\\
$g = 2x_1x_5x_6 + 2x_2x_6^2 - 2x_3x_5^2 + x_4x_9^2 - 2x_5x_8x_9 + 2x_6x_7x_9$\\$R(L)^L = k(x_4, f, g)$, which is rational over $k$.\\
\ \\
{\bf (4) \boldmath{$L_{9,4}$}}
\begin{eqnarray*}
M := \left[\begin{array}{ccccccccc}
0 &2x_2 &-2x_3 &0 &2x_5 &0 &-2x_7 &x_8 &-x_9\\
-2x_2 &0 &x_1 &0 &0 &2x_5 &x_6 &0 &x_8\\
2x_3 &-x_1 &0 &0 &x_6 &2x_7 &0 &x_9 &0\\
0 &0 &0 &0 &0 &0 &0 &0 &0\\
-2x_5 &0 &-x_6 &0 &0 &0 &0 &0 &0\\
0 &-2x_5 &-2x_7 &0 &0 &0 &0 &0 &0\\
2x_7 &-x_6 &0 &0 &0 &0 &0 &0 &0\\
-x_8 &0 &-x_9 &0 &0 &0 &0 &0 &x_4\\
x_9 &-x_8 &0 &0 &0 &0 &0 &-x_4 &0
\end{array}\right]
\end{eqnarray*}
$i(L) = 3$, $c(L) = 6$, $p_L = 1$, $F(L) = L$ (so $L$ is quasi quadratic with no CP's), $[L,L] = L$, $L_\lambda = L$,\\
$Y(L) = k[x_4, f, g]$ which is polynomial over $k$, where $f = 4x_5x_7 - x_6^2$,\\
$g = x_1x_4x_6 + 2x_2x_4x_7 - 2x_3x_4x_5 - x_5 x_9^2 + x_6x_8x_9 - x_7x_8^2$\\
$R(L)^L = k(x_4, f, g)$ which is rational over $k$.\\
\ \\
{\bf (5) \boldmath{$L_{9,5}$}}
\begin{eqnarray*}
M := \left[\begin{array}{ccccccccc}
0 &2x_2 &-2x_3 &x_4 &-x_5 &x_6 &-x_7 &0 &0\\
-2x_2 &0 &x_1 &0 &x_4 &0 &x_6 &0 &0\\
2x_3 &-x_1 &0 &x_5 &0 &x_7 &0 &0 &0\\
-x_4 &0 &-x_5 &0 &x_8 &0 &0 &0 &0\\
x_5 &-x_4 &0 &-x_8 &0 &0 &0 &0 &0\\
-x_6 &0 &-x_7 &0 &0 &0 &x_9 &0 &0\\
x_7 &-x_6 &0 &0 &0 &-x_9 &0 &0 &0\\
0 &0 &0 &0 &0 &0 &0 &0 &0\\
0 &0 &0 &0 &0 &0 &0 &0 &0
\end{array}\right]
\end{eqnarray*}
$i(L) = 3$, $c(L) = 6$, $p_L = 1$, $F(L) = L$ (so $L$ is quadi quadratic with no CP's, $[L,L] = L$, $L_\Lambda = L$,\\
$Y(L) = k[x_8, x_9, f]$, which is polynomial over $k$, where\\
$f = 2x_1x_4x_5x_9 + 2x_1x_6x_7x_8 + x_1^2x_8x_9 + 4x_2x_3x_8x_9 + 2x_2x_7^2x_8 + 2x_2x_5^2 x_9 - 2x_3x_4^2x_9 - 2x_3x_6^2x_8 + 2x_4x_5x_6x_7 - x_4^2x_7^2 - x_5^2x_6^2$\\
$R(L)^L = k(x_8, x_9, f)$ which is rational over $k$.\\
\ \\
{\bf (6) \boldmath{$L_{9,6}$}}
\begin{eqnarray*}
M := \left[\begin{array}{ccccccccc}
0 &2x_2 &-2x_3 &0 &0 &x_6 &-x_7 &x_8 &-x_9\\
-2x_2 &0 &x_1 &0 &0 &0 &x_6 &0 &x_8\\
2x_3 &-x_1 &0 &0 &0 &x_7 &0 &x_9 &0\\
0 &0 &0 &0 &0 &0 &0 &0 &0\\
0 &0 &0 &0 &0 &0 &0 &0 &0\\
-x_6 &0 &-x_7 &0 &0 &0 &0 &0 &x_4\\
x_7 &-x_6 &0 &0 &0 &0 &0 &-x_4 &0\\
-x_8 &0 &-x_9 &0 &0 &0 &x_4 &0 &x_5\\
x_9 &-x_8 &0 &0 &0 &-x_4 &0 &-x_5 &0
\end{array}\right]
\end{eqnarray*}
$i(L) = 3$, $c(L) = 6$, $p_L = 1$, $F(L) = L$ (so $L$ is quasi quadratic with no CP's), $[L,L] = L$, $L_\Lambda = L$,\\
$Y(L) = k[x_4, x_5, f]$, which is polynomial over $k$, where\\
$f = x_1^2x_4^2 + 2x_1x_4x_6x_9 + 2x_1x_4x_7x_8 - 2x_1x_5x_6x_7 + 4x_2x_4x_7x_9 - 2x_2x_5x_7^2 + 4x_2x_3x_4^2 - 4x_3x_4x_6x_8 + 2x_3x_5x_6^2 + x_6^2x_9^2 - 2x_6x_7x_8x_9 + x_7^2x_8^2$\\
$R(L)^L = k(x_4, x_5, f)$ which is rational over $k$.\\
\ \\
{\bf (7) \boldmath{$L_{9,7}$}}
\begin{eqnarray*}
M := \left[\begin{array}{ccccccccc}
0 &2x_2 &-2x_3 &2x_4 &0 &-2x_6 &2x_7 &0 &-2x_9\\
-2x_2 &0 &x_1 &0 &x_4 &2x_5 &0 &2x_7 &x_8\\
2x_3 &-x_1 &0 &2x_5 &x_6 &0 &x_8 &2x_9 &0\\
-2x_4 &0 &-2x_5 &0 &0 &0 &0 &0 &0\\
0 &-x_4 &-x_6 &0 &0 &0 &0 &0 &0\\
2x_6 &-2x_5 &0 &0 &0 &0 &0 &0 &0\\
-2x_7 &0 &-x_8 &0 &0 &0 &0 &x_4 &x_5\\
0 &-2x_7 &-2x_9 &0 &0 &0 &-x_4 &0 &x_6\\
2x_9 &-x_8 &0 &0 &0 &0 &-x_5 &-x_6 &0
\end{array}\right]
\end{eqnarray*}
$i(L) = 3$, $c(L) = 6$, $p_L = 1$, $F(L) = L$ (so $L$ is quasi quadratic with no CP's), $[L,L] = L$, $L_\lambda = L$,\\
$Y(L) = k[f,g,h]$, which is polynomial over $k$, where\\
$f = x_4x_6 - x_5^2$, $g = x_4x_9 - x_5x_8 + x_6x_7$\\
$h = 2x_1x_5 + 2x_2x_6 - 2x_3x_4 + 4x_7x_9 - x_8^2$\\
$R(L)^L = k(f,g,h)$, which is rational over $k$.\\
\ \\
In the following the radical $W$ of $L$ will be commutative and we denote by $W_m$ the standard $(m+1)$-dimensional irreducible $sl(2)$-module.\\
\ \\
{\bf (8) \boldmath{$L_{9,8} = sl(2) \oplus W$} with \boldmath{$W = W_1 \oplus W_1 \oplus W_1$}}
\begin{eqnarray*}
M := \left[\begin{array}{ccccccccc}
0 &2x_2 &-2x_3 &x_4 &-x_5 &x_6 &-x_7 &x_8 &-x_9\\
-2x_2 &0 &x_1 &0 &x_4 &0 &x_6 &0 &x_8\\
2x_3 &-x_1 &0 &x_5 &0 &x_7 &0 &x_9 &0\\
-x_4 &0 &-x_5 &0 &0 &0 &0 &0 &0\\
x_5 &-x_4 &0 &0 &0 &0 &0 &0 &0\\
-x_6 &0 &-x_7 &0 &0 &0 &0 &0 &0\\
x_7 &-x_6 &0 &0 &0 &0 &0 &0 &0\\
-x_8 &0 &-x_9 &0 &0 &0 &0 &0 &0\\
x_9 &-x_8 &0 &0 &0 &0 &0 &0 &0
\end{array}\right]
\end{eqnarray*}
$i(L) = 3$ ($=\dim W - \dim sl(2)$, so Proposition 1.10 is applicable.  In particular, $W$ is CPI).\\
$c(L) = 6$, $p_L = 1$, $F(L) = W$, $[L,L] = L$, $L_\Lambda = L$.\\
$Y(L) = k[f,g,h]$, which is polynomial over $k$, where
\begin{center}
$f = x_4x_7 - x_5x_6$, \ \ \ $g = x_4x_9 - x_5x_8$, \ \ \ $h = x_6x_9 - x_7x_8$
\end{center}
$R(L)^L = k(f,g,h)$, which is rational over $k$.\\
\ \\
{\bf (9) \boldmath{$L_{9,9} = sl(2) \oplus W$} with \boldmath{$W = W_2 \oplus W_2$}}\\
This is the smallest counterexample to the Gelfand Kirillov conjecture [1,2]. See also [45, Example 60].
\begin{eqnarray*}
M := \left[\begin{array}{ccccccccc}
0 &2x_2 &-2x_3 &2x_4 &0 &-2x_6 &2x_7 &0 &-2x_9\\
-2x_2 &0 &x_1 &0 &2x_4 &x_5 &0 &2x_7 &x_8\\
2x_3 &-x_1 &0 &x_5 &2x_6 &0 &x_8 &2x_9 &0\\
-2x_4 &0 &-x_5 &0 &0 &0 &0 &0 &0\\
0 &-2x_4 &-2x_6 &0 &0 &0 &0 &0 &0\\
2x_6 &-x_5 &0 &0 &0 &0 &0 &0 &0\\
-2x_7 &0 &-x_8 &0 &0 &0 &0 &0 &0\\
0 &-2x_7 &-2x_9 &0 &0 &0 &0 &0 &0\\
2x_9 &-x_8 &0 &0 &0 &0 &0 &0 &0
\end{array}\right]
\end{eqnarray*}
$i(L) = 3$ ($=\dim W - \dim sl(2)$, so by Proposition 1.10 $W$ is a CPI), $c(L) = 6$, $p_L = 1$, $F(L) = W$, $[L,L]= L$, $L_\Lambda = L$.\\
$Y(L) = k [f,g,h]$, which is polynomial over $k$, where
\begin{center}
$f = 4x_4x_6 - x_5^2$, \ \ \ $g = 4x_7x_9 - x_8^2$, \ \ \ $h = 2x_4x_9 + 2x_6x_7 - x_5x_8$
\end{center}
$R(L)^L = k(f,g,h)$, which is rational over $k$.\\
\ \\
{\bf (2b) \boldmath{$L$} is not coregular}\\
{\bf (10) \boldmath{$L_{9,10} = sl(2) \oplus W$} with \boldmath{$W = W_3 \oplus W_1$}}
\begin{eqnarray*}
M := \left[\begin{array}{ccccccccc}
0 &2x_2 &-2x_3 &3x_4 &x_5 &-x_6 &-3x_7 &x_8 &-x_9\\
-2x_2 &0 &x_1 &0 &3x_4 &2x_5 &x_6 &0 &x_8\\
2x_3 &-x_1 &0 &x_5 &2x_6 &3x_7 &0 &x_9 &0\\
-3x_4 &0 &-x_5 &0 &0 &0 &0 &0 &0\\
-x_5 &-3x_4 &-2x_6 &0 &0 &0 &0 &0 &0\\
x_6 &-2x_5 &-3x_7 &0 &0 &0 &0 &0 &0\\
3x_7 &-x_6 &0 &0 &0 &0 &0 &0 &0\\
-x_8 &0 &-x_9 &0 &0 &0 &0 &0 &0\\
x_9 &-x_8 &0 &0 &0 &0 &0 &0 &0
\end{array}\right]
\end{eqnarray*}
 $i(L) = 3$, $c(L) = 6$, $p_L = 1$, $F(L) = W$, $[L,L] = L$, $L_\Lambda = L$.\\
Clearly $i(L) = \dim W -3$ and so by Proposition 3.6 $W$ is a CP-ideal of $L$ and $L$ satisfies Dixmier's 4th problem.\\
Let us now demonstrate that $L$ is not coregular.  Suppose on the contrary that $Y(L)$ is freely generated by some homogeneous $f_1, f_2, f_3$.  Note that we have the following equality:
$$3i(L) + 2\deg p_L = 9 = \dim L + 2\dim Z(L)$$
($\deg p_L = 0$ and $\dim Z(L) = 0$).  This implies that $\deg f_i \leq 2$, $i = 1,2,3$ by Corollary 1.8.  However using MAPLE, it is easy to see that there are no such invariants. In fact 4 is the smallest degree of a nontrivial invariant. Contradiction.\\
\ \\
{\bf (11) \boldmath{$L_{9,11} = sl(2) \oplus W_5$}}
\begin{eqnarray*}
M := \left[\begin{array}{ccccccccc}
0 &2x_2 &-2x_3 &5x_4 &3x_5 &x_6 &-x_7 &-3x_8 &-5x_9\\
-2x_2 &0 &x_1 &0 &5x_4 &4x_5 &3x_6 &2x_7 &x_8\\
2x_3 &-x_1 &0 &x_5 &2x_6 &3x_7 &4x_8 &5x_9 &0\\
-5x_4 &0 &-x_5 &0 &0 &0 &0 &0 &0\\
-3x_5 &-5x_4 &-2x_6 &0 &0 &0 &0 &0 &0\\
-x_6 &-4x_5 &-3x_7 &0 &0 &0 &0 &0 &0\\
x_7 &-3x_6 &-4x_8 &0 &0 &0 &0 &0 &0\\
3x_8 &-2x_7 &-5x_9 &0 &0 &0 &0 &0 &0\\
5x_9 &-x_8 &0 &0 &0 &0 &0 &0 &0
\end{array}\right]
\end{eqnarray*}
$i(L) = 3$, $c(L) = 6$, $p_L = 1$, $F(L) = W_5$, $W_5$ is a CPI, $[L,L] = L$, $L_\Lambda = L$.\\
$L$ satisfies Dixmier's 4th problem by (9) of Theorem 2.1 (or by Proposition 3.6 because $i(L) = \dim W_5 - 3$).  The fact that $L_{9,11}$ is not coregular can be shown in the same fashion as $L_{9,10}$.  Moreover,
\begin{center}
$Y(L_{9,11}) = S(W_5)^{sl(2)} = S(W_5)^{SL(2)}$ by Proposition 1.10
\end{center}
Now suppose $k = \Bbb C$.  Then $W_5$ may be considered as the vector space of binary forms of degree 5 with complex coefficients (the quintics) on which $SL(2,\Bbb C)$ acts.
The algebra of invariants $\Bbb C[W_5]^{SL(2,\Bbb C)}$ which is isomorphic to $S(W_5)^{sl(2,\Bbb C)}$, has been studied already in the 19th century by Sylvester, among others.  
At first 3 algebraically independent
invariants $I_4$, $I_8$, $I_{12}$ were found of degrees 4, 8, 12.  In 1854 Hermite discovered an invariant $I_{18}$ of degree 18 and he showed that
$$\Bbb C[W_5]^{SL(2,\Bbb C)} = \Bbb C[I_4, I_8, I_{12}, I_{18}]$$
with the following relation:
$$16I_{18}^2 = I_4I_8^4 + 8I_8^3I_{12} - 2I_4^2I_8^2I_{12} - 72I_4I_8I_{12}^2 - 432I_{12}^3 + I_4^3 I_{12}^2$$
In particular, $\Bbb C[W_5]^{SL(2,\Bbb C)}$ is not polynomial.\\
The explicit forms of $I_4$, $I_8$, $I_{12}$, $I_{18}$ were given in papers by Cayley.  $I_{18}$ has 848 monomials with very large coefficients ! See [12, p.41]. 
\hfill $\square$\\
\ \\
{\bf Remark 4.7}  In higher dimensions Dixmier's 4th problem will become more challenging.  For instance, it may occur that for certain Lie algebras the problem is precisely equivalent to the rationality problem in invariant theory, which is still open.\\
\ \\
{\bf Acknowledgments}\\
\\
We would like to thank Michel Van den Bergh for some helpful discussions and for providing Proposition 3.1.  We are also grateful to Vladimir Popov for accurately describing the relevant results (and their references) of Katsylo and Bogomolov.\\
Finally we wish to thank Viviane Mebis for the excellent typing of the manuscript.

\renewcommand\bibname{\normalsize{References}}


\begin{thebibliography}{AA}
\bibitem{1} J. Alev, A. Ooms and M. Van den Bergh, A class of counterexamples to the Gelfand-Kirillov conjecture, Trans. Amer. Math. Soc. 348 (1996), 1709-1716.
\bibitem{2} J. Alev, A. Ooms and M. Van den Bergh, The Gelfand-Kirillov conjecture for Lie algebras of dimension at most eight, J. Algebra 227 (2000), 549-581. Corrigendum J. Algebra 230 (2000), 749.
\bibitem{3} P. Bernat, Sur le corps enveloppant d'une alg\`ebre de Lie r\'esoluble, C.R. Acad. Sci. Paris 258 (1964), 2713-2715.
\bibitem{4} F.A. Bogomolov, P.I. Katsylo, Rationality of some quotient varieties, Math. USSR Sb 54 (1986) 571-576.
\bibitem{5} W. Borho, P. Gabriel, R. Rentschler, Primideale in Einh\"{u}llenden aufl\"{o}sbarer Lie-Algebren, Lecture Notes in Math., vol. 357, Springer-Verlag, Berlin, 1973.
\bibitem{6} N. Boubaki, ``Groupes et alg\`ebres de Lie'', Chap. I, Hermann, Paris, 1971.
\bibitem{7} V. Coll, A. Giaquinto, C. Magnant, Meander graphs and Frobenius seaweed algebras, J. Gen. Lie Theory Appl. 5 (2011), 1-7. 
\bibitem{8} V. Coll, C. Magnant, H. Wang, The signature of a meander, arXiv:1206.2705v2 [math.QA], 2014.
\bibitem{9} B. Csikos, L. Verhoczki, Classification of  Frobenius Lie algebras of dimension $\leq$ 6, Publ. Math. Debrecen, 70 (2007) 427-451.
\bibitem{10} L. Delvaux, E. Nauwelaerts, A.I. Ooms, On the semicenter of a universal enveloping algebra, J. Algebra 94 (1985), 324-346.
\bibitem{11} J. Dixmier, Sur le centre de l' alg\`ebre enveloppante d'une alg\`ebre de Lie, C. R. Acad. Sci. Paris A 265 (1967), 408-410.
\bibitem{12} J. Dixmier, Quelques aspects de la th\'eorie des invariants, Gaz. Math., Soc. Math. Fr., 43 (1990), 39-64.
\bibitem
{13} J. Dixmier, ''Enveloping Algebras'', Grad. Stud. Math., vol 11, Amer. Math. Soc., Providence, RI, 1996.
\bibitem{14} J. Dixmier, M. Duflo, M. Vergne, Sur la repr\'esentation coadjointe d' une alg\`ebre de Lie, Compositio Math. 29 (1974), 309-323.
\bibitem{15} M. Duflo, Constructions of primitive ideals in an enveloping algebra, Publ. of 1971 Summer School in Math., edited by I.M. Gelfand, Bolyai-Janos Math. Soc., Budapest.
\bibitem{16} M. Duflo, Sur les extensions des repr\'esentations irr\'eductibles des groupes de Lie nilpotents, Ann. Scient Ec. Norm. Sup. 5 (1972), 71-120.
\bibitem{17} M. Duflo, Op\'erateurs diff\'erentielles sur un groupe de Lie, Ann. Scient Ec. Norm. Sup. 10 (1977), 265-288.
\bibitem{18} M. Duflo, R.W.T. Yu, On compositions associated to Frobenius parabolic and seaweed subalgebras of sln(k), Journal of Lie Theory 25 (2015), 1191-1213. 
\bibitem{19} A. G. Elashvili, Frobenius Lie algebras, Funct. Anal. i Prilozhen 16 (1982), 94-95.
\bibitem{20} A. G. Elashvili, Frobenius Lie algebras II, Trudy Razmadze Math. Institute (Tbilisi), 77 (1985), 127-137.
\bibitem{21} A. G. Elashvili, On the index of orispherical subalgebras of semisimple Lie algebras, Trudy Razmadze Math. Institute (Tbilisi) 77 (1985), 116-126.
\bibitem{22} A. G. Elashvili, A. I. Ooms, On commutative polarizations,  J. Algebra 264  (2003), 129-154.
\bibitem{23} F. Fauquant-Millet, Sur la polynomiality de certaines alg\`ebres d'invariants d'alg\`ebres de Lie, M\'emoire d'Habilitation \'a Diriger des Recherches (2014), https://tel.archives-ouvertes.fr/tel-00994655.
\bibitem{24} F. Fauquant-Millet, A. Joseph, Semi-centre de l'alg\`ebre enveloppante d' une sous-alg\`ebre parabolique d' une alg\`ebre de Lie semi-simple, Ann. Sci. \'Ecole Norm. Sup. 38 (2005), 155-191.
\bibitem{25} F. Fauquant-Millet, A. Joseph, La somme des faux degr\'es - un myst\`ere en th\'eorie des invariants, Adv. Math. 217 (2008), 1476-1520.
\bibitem{26} F. Fauquant-Millet, P. Lamprou, Polynomiality for the Poission centre of truncated maximal parabolic subalgebras arXiv:1701.02238.
\bibitem{27} I.M. Gelfand, A.A. Kirillov, Sur les corps li\'es aux alg\`ebres enveloppantes des alg\`ebres de Lie, Inst. Hautes Etudes Sci. Publ. Math. 31 (1966), 5-19.
\bibitem{28} A. Joseph, On semi-invariants and index for biparabolic (seaweed) algebras. I, J. Algebra 305 (2006) 487-515.
\bibitem{29} A. Joseph, On semi-invariants and index for biparabolic (seaweed) algebras. II, J. Algebra 312 (2007) 158-193.
\bibitem{30} A. Joseph, The hidden semi-invariants generators of an almost-Frobenius biparabolic, Transform. Groups, 19 (2014), 735-778.
\bibitem{31} A. Joseph, Slices for biparabolic coadjoint actions in type A, J. Algebra 319 (2008), 5060-5100.
\bibitem{32} A. Joseph, An algebraic slice in the coadjoint space of the Borel and the Coxeter element, Advances in Mathematics 227 (2011), 522-585.
\bibitem{33} A. Joseph, D. Shafrir, Polynomiality of invariants, unimodularity and adapted pairs. Transform. Groups, 15 (2010), 851-882.
\bibitem{34} P.I. Katsylo, Rationality of fields of invariants of reducible representations of $SL(2)$, Mosc. Univ. Math. Bull 39 (1984), no. 5.
\bibitem{35} L. Le Bruyn, A.I. Ooms, The semicenter of an enveloping algebra is factorial, Proc. Amer. Mat. Soc. 93 (1985), 397-400.
\bibitem{36}C. Moeglin, Factorialit\'e dans les alg\`ebres enveloppantes, C.R. Acad. Sci. Paris S\'erie A 282 (1976),1269-1272.
\bibitem{37}T. Moons, A.I. Ooms, On the Jordan kernel of a universal enveloping algebra, J. Algebra 122 (1989), 211-231.
\bibitem{38} A.I. Ooms, On Lie algebras having a primitive universal enveloping algebra, J. Algebra, 32 (1974), 488-500.  
\bibitem{39} A.I. Ooms, On Lie algebras with primitive envelopes, supplements, Proc. Amer. Math. Soc. 58 (1976), 67-72.
\bibitem{40} A.I. Ooms, On Frobenius Lie algebras, Comm. Algebra 8 (1980), 13-52.
\bibitem{41} A.I. Ooms, On certain maximal subfields in the quotient division ring of an enveloping algebra, J. Algebra, 230 (2000), 694-712.
\bibitem{42} A.I. Ooms, The Frobenius semiradical of a Lie algebra, J. Algebra 273 (2004), 274-287.
\bibitem{43} A.I. Ooms, Computing invariants and semi-invariants by means of Frobenius Lie algebras,  J. Algebra 321 (2009), 1293-1312.
\bibitem{44} A.I. Ooms, The Poisson center and polynomial, maximal Poisson commutative subalgebras, especially for nilpotent Lie algebras of dimension at most seven, J. Algebra 365 (2012), 83-113.
\bibitem{45} A.I. Ooms, The polynomiality of the Poisson center and semi-center of a Lie algebra and Dixmier's fourth problem, J. Algebra 477 (2017) 95-146.
\bibitem{46} A.I. Ooms, The polynomiality of the Poisson center and semi-center of a Lie algebra and Dixmier's fourth problem, arXiv:1605.04200.
\bibitem{47} A.I. Ooms, M. Van den Bergh, A degree inequality for Lie algebras with a regular Poisson semicenter, J. Algebra 323 (2010), 305-322.
\bibitem{48} D.I. Panyushev, An extension of Ra\"{i}s theorem and seaweed subalgebras of simple Lie algebras, Ann. Inst. Fourier (Grenoble) 55 (3) (2005), 693-715.
\bibitem{49} D.I. Panyushev, O.S. Yakimova, On seaweed subalgebras and meander graphs in type C, arXiv:math.RT/1601.00305v1 (2016).
\bibitem{50} R. Rentschler, M. Vergne, Sur le semi-centre du corps enveloppant d'une alg\`ebre de Lie, Ann. Sci. \'Ecole Norm. Sup. 6 (1973), 389-405.
\bibitem{51} D. Shafrir, Private communication.
\bibitem{52} P. Tauvel, R. T. W. Yu, Lie Algebras and Algebraic Groups, Springer Monographs in Mathematics, Springer-Verlag, Berlin, 2005.
\bibitem{53} P. Tauvel, R. T. W. Yu, Affine slice for the coadjoint action of a class of biparabolic subalgebras of a semisimple Lie algebra, Algebr. Represent. Theor. 16 (2013), 859-872.
\bibitem{54} M. Vergne, La structure de Poisson sur l' alg\`ebre sym\'etrique d'une alg\`ebre de Lie nilpotente, Bull. Soc. Math. France, 100 (1972), 301-335.
\bibitem{55} O. Yakimova, A counterexample to Premet's and Joseph's conjectures, Bull. London Math. Soc., 39 (2007), 749-754.
\bibitem{56} A.I. Ooms, The maximal abelian dimension of a Lie algebra, Rentschler's property and Milovanov's conjecture, Algebr. Represent. Theor. (2019).  https://doi.org/10.1007/s10468-019-09877-5.
\bibitem{57} F. Fauquant-Millet, Weierstrass sections for some truncated parabolic subalgebras, arXiv: 1907.11755 [math.RT].
\end{thebibliography}
\end{document}